\documentclass[11pt, a4paper]{article}
\author{\textit{I. V.} \textit{Protasov,} \textit{A.} \textit{Tsvietkova} }
\title{\textbf{Decomposition} \textbf{of Cellular Balleans}}
\date{}
\usepackage{formula}
\usepackage{latexsym}
\usepackage{amsmath}
\usepackage{amsfonts}
\usepackage{eucal}
\begin{document}
\maketitle
 \footnotesize
 \begin{description} \item \qquad \textbf{Abstract.}  A ballean is a set endowed with some family of its subsets which are called the balls. We
 postulate the properties of the family of balls in such a way
 that the balleans can be considered as the asymptotic
 counterparts of the uniform topological spaces. The isomorphisms
 in the category of balleans are called asymorphisms. Every metric
 space can be considered as a ballean. The ultrametric spaces are
 prototypes for the cellular balleans. We prove some general
 theorem about decomposition of a homogeneous cellular ballean in a
 direct product of a pointed family of sets. Applying this
 theorem we show that the balleans of two uncountable groups of
 the same regular cardinality are asymorphic.\\ 
\end{description} \normalsize \

A \textit{ball structure} is a triple $\mathcal{B}=(X,P,B)$ where
$X,P$ are non-empty sets, and for all $ x\in X $ and $ \alpha\in P
$ , $ B(x,\alpha) $ is a subset of $X$ which is called a
\textit{ball of radius} $\alpha$ around $x$. It is supposed that
$x \in B(x,\alpha)$ for all $x \in X$, $ \alpha \in P$. The set
$X$ is called the \textit{support} of $\mathcal{B}$, $P$ is called
the \textit{set of radii}. \

Given any $x \in X , A\subseteq X, \alpha \in P$, we put
\begin{eqnarray} \nonumber B^{*}(x,\alpha)=\{y \in X : x \in B(y,\alpha)\}, \\
\ \displaystyle B(A,\alpha)=\bigcup_{a \in A}B(a,\alpha),\\
B^*(A,\alpha)=\bigcup_{a \in A}B^*(a,\alpha). \
\end{eqnarray}

A ball structure $\mathcal{B}=(X,P,B)$ is called \textit{a
ballean} (or a \textit{a coarse structure}) if
\begin{itemize}
\item $\forall \alpha, \beta \in P \ \exists \alpha', \beta' \in
P$ such that $\forall x \in X$
\begin{eqnarray}
 \nonumber B(x,\alpha)\subseteq B^{*}(x,\alpha'), \ B^{*}(x,\beta)\subseteq
 B(x,\beta');
\end{eqnarray}
\item $\forall \alpha,\beta \in P\ \exists \gamma \in P$ such that
$\forall x \in X$
\begin{eqnarray}
\nonumber B(B(x,\alpha),\beta)\subseteq B(x,\gamma);
\end{eqnarray}
\end{itemize} \qquad

Let $\mathcal{B}_{1}=(X_{1},P_{1},B_{1})$ and
$\mathcal{B}_{2}=(X_{2},P_{2},B_{2})$ be balleans.\

A mapping $f:X_{1}\rightarrow X_{2}$ is called a
\textit{$\prec$-mapping} if $\forall \alpha \in P_{1}\ \exists
\beta \in P_{2}$ such that:
\begin{eqnarray}
\nonumber f(B_{1}(x,\alpha))\subseteq B_{2}(f(x),\beta).
\end{eqnarray}\ \qquad
A bijection $f:X_{1}\rightarrow X_{2}$ is called an
\textit{asymorphism} between $\mathcal{B}_{1}$ and
$\mathcal{B}_{2}$ if $f$ and $f^{-1}$ are $\prec$-mappings. In
this case $\mathcal{B}_{1}$ and $\mathcal{B}_{2}$ are called
\textit{asymorphic}.\ \

If $X_1=X_2$ and the identity mapping id: $X_{1}\rightarrow X_{2}$
is an asymorphism, we identify $\mathcal{B}_{1}$ and
$\mathcal{B}_{2}$ and write $\mathcal{B}_{1}=\mathcal{B}_{2}$.

For motivation to study balleans, see [1], [2], [3], [4].
\

Every metric space $(X,d)$ determines the \textit{metric ballean}
$\mathcal{B}(X,d)$\

$=(X,\mathbb{R^{+}},B_{d})$, where
$\mathbb{R^{+}}$ is the set of non-negative real numbers, \
\begin{center}
$B_{d}(x,r)=\{y\in X: d(x,y)\leq r \}$.\
\end{center}
A ballean $\mathcal{B}$ is called \textit{metrizable} if
$\mathcal{B}$ is asymorphic to $\mathcal{B}(X,d)$ for some metric
ballean. By [3,Theorem 2.1], a ballean $\mathcal{B}$ is metrizable
if and only if $\mathcal{B}$ is connected and the cofinality
cf($\mathcal{B})\leq\aleph_0$. A ballean $\mathcal{B}=(X,P,B)$ is
\textit{connected} if, for any $x,y\in X$, there exists $\alpha\in
P$ such that $y\in B(x, \alpha)$. To define cf($\mathcal{B})$, we
use the natural preordering on $P$: $\alpha\leq\beta$ if and only
if $B(x,\alpha)\subseteq B(x, \beta)$ for every $x\in X$. A subset
$P'$ is \textit{cofinal} in $P$ if, for every $\alpha\in P$, there
exists $\alpha' \in P'$ such that $\alpha \leq\ \alpha'$, so
cf($\mathcal{B})$ is the minimal cardinality of cofinal subsets of
$P$. \

Given an arbitrary ballean $\mathcal{B}=(X,P,B)$, $x,y \in X$ and
$\alpha \in P$, we say that
 $x,y$ are \textit{$\alpha$-path connected} if there exists a finite
sequence $x_{0}, x_{1}, ..., x_{n}$, $x_{0}=x$, $x_{n}=y$  such
that $x_{i+1} \in B(x_{i},\alpha)$, for every $i \in \{0,1,...,n-1
\}$. For any $x \in X$ and $\alpha \in P$, we put \
\\
\begin{center} $B^{\Box}(x,\alpha)=\{y \in X : x,y$ are $\alpha$-path
connected$\}$ \end{center}

 The ballean $\mathcal{B}^{\Box}=(X,P,B^{\Box})$ is called the
 \textit{cellularization} of $\mathcal{B}$. A ballean
 $\mathcal{B}$ is called \textit{cellular} if
 $\mathcal{B}^{\Box}=\mathcal{B}$. For characterizations of
 cellular balleans see [3, Chapter 3].

\

 \textbf{Example 1}. A metric $d$ on a set $X$ is called an
 \textit{ultrametric} if
 \
 \begin{center} $d(x,y)\leq$ max$\{d(x,z), d(y,z)\}$ \end{center}
 for all $x,y,z \in X$. If (X,d) is an ultrametric space then the
 ballean $\mathcal{B}(X,d)$ is cellular. Moreover, by [3, Theorem
 3.1], a ballean $\mathcal{B}$ is metrizable and cellular if
 and only if $\mathcal{B}$ is asymorphic to the metric ballean $\mathcal{B}(X,d)$
of some ultrametric space $(X,d)$.

\

\textbf{Example 2}. Let $G$ be an infinite group with the identity
$e$, $\kappa$ be an infinite cardinal such that $\kappa \leq |G|$,
$\mathcal{F}(G,\kappa)=\{A\subseteq G:e\in A, |A|<\kappa\}$. Given
any $g\in G$ and $A\in F(G,\kappa)$, we put $B(g,A)=gA$ and get
the ballean $\mathcal{B}(G,\kappa)=(G, F(G, \kappa), B)$. In the
case $\kappa =|G|$, we write $\mathcal{B}(G)$ instead of
$\mathcal{B}(G,\kappa)$. A ballean $\mathcal{B}(G,\kappa)$ is
cellular if and only if either $\kappa>\aleph_0$ or
$\kappa=\aleph_0$ and $G$ is locally finite (i.e. every finite
subset of $G$ is contained in some finite subgroup).

\

\textbf{Example 3}. A family of subsets of a group $G$ is called a
\textit{Boolean group ideal} if
\begin{itemize}
\item $A,B \in \Im \Rightarrow A\cup B \in\Im$; \item $A \in \Im,
A'\subset A \Rightarrow A' \in\Im$; \item $A,B \in \Im \Rightarrow
AB \in\Im, A^{-1}\in \Im$; \item $F \in \Im$ for every finite
subset $F$ of $G$.
\end{itemize} \qquad
\ Every Boolean group ideal $\Im$ on $G$ determines the ballean
$\mathcal{B}(G,\Im)=(G,\Im,B)$, where $B(g,A)=gA$ for all $g \in
G, A\in \Im$. The balleans on groups determined by the Boolean
group ideals can be considered (see [3, Chapter 6]) as the
asymptotic counterparts of the group topologies. A ballean
$\mathcal{B}(G,\Im)$ is cellular if and only if $\Im$ has a base
consisting of the subgroups of $G$.

\ 

A connected ballean $\mathcal{B}=(X,P,B)$ is called ordinal if
there exists a cofinal well-ordered (by $\leq$) subset of $P$.
Clearly, every metrizable ballean is ordinal.

 \

\textbf{Theorem 1}.\textit{ Let $\mathcal{B}=(X,P,B)$
 be an ordinal ballean. Then $\mathcal{B}$ is either metrizable or
 cellular.}

\

\textit{Proof}. If cf($\mathcal{B})\leq\aleph_0$ then $\mathcal{B}$
is metrizable by theorem 2.1 from [3]. Assume that
cf($\mathcal{B})>\aleph_0$. Given an arbitrary $\alpha\in P$, we
choose inductively a sequence $(\alpha_n)_{n\in\omega}$ in $P$
such that $\alpha_0=\alpha$ and $B(B(x,\alpha_n),\alpha)\subseteq
B(x,\alpha_{n+1})$ for every $x\in X$. Since
cf($\mathcal{B})>\aleph_0$, we can pick $\beta\in P$ such that
$\beta\geq\alpha_n$ for every $n\in\omega$. Then $B^{\Box}(x,
\alpha)\subseteq B(x, \beta)$ for every $x\in X$, so
$\mathcal{B}^{\Box}=\mathcal{B}$.

\

Let $\gamma$ be an ordinal, $\{Z_{\lambda}:\lambda<\gamma\}$ be a
family of non-empty sets. For every $\lambda<\gamma$  we fix some
element $e_{\lambda}\in Z_{\lambda}$ and say that the family
$\{(Z_{\lambda}, e_{\lambda}):\lambda<\gamma\}$ is
\textit{pointed}. A \textit{direct product}
$Z=\otimes_{\lambda<\gamma}(Z_{\lambda}, e_{\lambda})$ is the set
of all functions $f:\{\lambda:\lambda<\gamma\}\rightarrow
\cup_{\lambda<\gamma}Z_{\lambda}$ such that $f(\lambda)\in
Z_{\lambda}$ and $f(\lambda)=e_{\lambda}$ for all but finitely
many $\lambda<\gamma$. We consider the ball structure
$\mathcal{B}(Z)=(Z,\{\lambda:\lambda<\gamma\},B)$, where
$B(f,\lambda)=\{g\in Z: f(\lambda')=g(\lambda')$ for all
$\lambda'\geq\lambda\}$ It is easy to verify that $\mathcal{B}(Z)$
is a cellular ballean. \

We say that a ballean $\mathcal{B}$ is \textit{decomposable in a
direct product} if $\mathcal{B}$ is asymorphic to $\mathcal{B}(Z)$
for some direct product $Z$.

\

 \textbf{Theorem 2}.\textit{ Let
$\gamma$ be a limit ordinal,
$\mathcal{B}=(Z,\{\lambda:\lambda<\gamma\},B)$ be a ballean such
that:
\begin{description}
\item[$(i)$] $B^{\Box}(x,\alpha)=B(x,\alpha)$ for all $x\in X,
\alpha \in P$; \item[$(ii)$] if $\alpha<\beta<\gamma$ then
$B(x,\alpha)\subset B(x,\beta)$ for each $x\in X$; \item[$(iii)$]
if $\beta$ is a limit ordinal and $\beta<\gamma$ then
$B(x,\beta)=\cup_{\alpha<\beta}B(x,\alpha)$ for each $x\in X$;
\item[$(iv)$] there exists a cardinal $\kappa_{0}$ such that
$B(x,0)=\kappa_{0}$ for each $x\in X$; \item[$(v)$]for every
$\alpha<\gamma$ there exists a cardinal $\kappa_{\alpha}$ such
that every ball of radius $\alpha+1$ is a disjoint union of
$\kappa_{\alpha}$-many balls of radius $\alpha$. \end{description}
\
Then
$\mathcal{B}$ is decomposable in a direct product.}

\

\textit{Proof}. We fix some set $Z_{0}$ of cardinality
$\kappa_{0}$ and define inductively a family of sets
$\{Z_{\alpha},\alpha<\gamma\}$. If $\alpha$ is a limit ordinal, we
take $Z_{\alpha}$ to be a singleton. If $\alpha=\beta+1$ we take a
set $Z_{\alpha}$ of cardinality $\kappa_{\beta}$. For every
$\alpha<\gamma$, we choose some element $e_{\alpha}\in
Z_{\alpha}$, put $Z=\otimes_{\lambda<\gamma}(Z_{\lambda},
e_{\lambda})$ and show that $\mathcal{B}$ is asymorphic to
$\mathcal{B}(Z)$. To this end we fix some element $x_{0}\in X$
and, for every $\alpha<\gamma$, define a mapping
$f_{\alpha}:B(x_{0},\alpha)\rightarrow\otimes_{\beta\leq\alpha}(Z_{\beta},e_{\beta})$
such that, for all $\beta<\alpha<\gamma,
f_{\alpha}\mid_{B(x_{0},\beta)}=f_{\beta}$ and the inductive limit
$f$ of the family $\{f_{\alpha}:\alpha<\gamma\}$ is an asymorphism
between $\mathcal{B}$ and $\mathcal{B}(Z)$. Here we identify
$\otimes_{\beta\leq\alpha}(Z_{\beta},e_{\beta})$ with the
corresponding subset of
$\otimes_{\beta<\gamma}(Z_{\beta},e_{\beta})$.

\ At the first step we fix some bijection
$f_{0}:B(x_{0},0)\rightarrow Z_{0}$ such that
$f_{0}(x_{0})=e_{0}$. Let us assume that, for some
$\alpha<\gamma$, we have defined the mappings
$\{f_{\beta}:\beta<\alpha\}$. If $\alpha$ is a limit ordinal, we
put $f_{\alpha}:
B(x_{0},\alpha)\rightarrow\otimes_{\beta<\alpha}(Z_{\beta},e_{\beta})$
to be an inductive limit of the family
$\{f_{\beta}:\beta<\alpha\}$. Since $Z_{\alpha}=\{e_{\alpha}\}$ we
can identify $\otimes_{\beta<\alpha}(Z_{\beta},e_{\beta})$ with
$\otimes_{\beta\leq\alpha}(Z_{\beta},e_{\beta})$, so $f_{\alpha}:
B(x_{0},\alpha)\rightarrow\otimes_{\beta\leq\alpha}(Z_{\beta},e_{\beta})$.
If $\alpha=\beta+1$, by cellularity of $\mathcal{B}$, there exists
a subset $Y\subseteq B(x_{0},\alpha), x_{0}\in Y$ such that
$B(x_{0},\alpha)$ is a disjoint union of the family
$\{B(y,\beta):y\in Y\}$. For every $y\in Y$, we can repeat the
inductive procedure of construction of $f_{\alpha}:
B(x_{0},\beta)\rightarrow \otimes_{\lambda \leq
\beta}(Z_{\lambda},e_{\lambda})$ to define a mapping $f'_{\beta,
y}B(y,\beta)\rightarrow\otimes_{\lambda\leq\beta}(Z_{\lambda},e_{\lambda})$.
Thus we fix some bijection $h: Y\rightarrow Z_{\alpha},
h(x_{0})=e_{\alpha}$ and put $f_{\beta,y}(x)=(f'_{\beta,y}(x),
h(y)), x\in B(y,\beta)$. At last, given any $x\in
B(x_{0},\alpha)$, we choose $y\in Y$ such that $x\in B(y,\beta)$
and put $f_{\alpha}(x)=f_{\beta,y}(x)$. By the construction of $f$
as an inductive limit of the family
$\{f_{\alpha}:\alpha<\gamma\}$, given any $x\in X$ and
$\alpha<\gamma$, we have $f(B(x,\alpha))=B(f(x),\alpha)$ so $f$ is
an asymorphism.

 \

 In the next two corollaries and Theorem 3 $\mathcal{B}(G)$ is a
ballean defined in Example 2.

\

\textbf{Corollary 1.} \textit{ Let $G$ be a countable locally
finite group. Then $\mathcal{B}(G)$ is decomposable in a direct
product of finite sets.}

\

\textit{Proof}. We write $G$ as a union $G=\cup_{n<\omega}G_{n}$
of an increasing chain of finite groups. Clearly, $\mathcal{B}(G)$
is asymorphic to the ballean $\mathcal{B}=(G,\omega,B)$ where
$B(g,n)=gG_{n}$. We put $\kappa_{0}=|G_{0}|$,
$\kappa_{n+1}=|G_{n+1}:G_{n}|$ and apply Theorem 2. 

\

\textbf{Corollary 2.} \textit{ Let $G$ be an uncountable group of
regular cardinality $\gamma$. Then $\mathcal{B}(G)$ is
decomposable in a direct
product.}

\

\textit{Proof}. We write $G$ as a union
$G=\cup_{\alpha<\gamma}G_{\alpha}$ of an increasing chain of
subgroups such that $|G_{0}|=\aleph_{0}$, $|G_{\alpha}|<\gamma$
and $G_{\alpha}=\cup_{\beta<\alpha}G_{\beta}$ for every limit
ordinal $\alpha$. Since $\gamma$ is regular, every subset
$F\subset G$, $|F|<|G|$ is contained in some subgroup
$G_{\alpha}$. It follows that $\mathcal{B}(G)$ is asymorphic to
the ballean $\mathcal{B}=(G,\gamma,B)$, where
$B(g,\alpha)=gG_{\alpha}$. Apply Theorem 2. 

 \

\textbf{Theorem 3}.\textit{ Let $G, H$ be two uncountable groups
of the same regular cardinality $\gamma$. Then $\mathcal{B}(G)$
and $\mathcal{B}(H)$ are asymorphic.}

\

 \textit{Proof}. We consider
two cases.\

\textit{Case 1:} $\gamma$ is a limit cardinal. We choose an
increasing family $\{G_{\alpha}:\alpha<\gamma\}$ of subgroups of
$G$ such that $G=\cup_{\alpha<\gamma}G_{\alpha}$,
$|G_{0}|=\aleph_{0}$, $|G_{\alpha+1}|=|G_{\alpha}|^{+}$ and
$G_{\beta}=\cup_{\alpha<\beta}G_{\alpha}$ for every limit ordinal
$\beta$. Put $\kappa_{\alpha}=|G_{\alpha}|^{+}, \alpha<\gamma$. By
Theorem 2, $\mathcal{B}(G)$ is asymorphic to $\mathcal{B}(Z)$
where the direct product $Z$ is defined by the family of cardinals
$\{\kappa_{\alpha}:\alpha<\gamma\}$. Since $H$ admits a filtration
$H=\cup_{\alpha<\gamma}H_{\alpha}$ with the same family
$\{\kappa_{\alpha}:\alpha<\gamma\}$ of parameters,
$\mathcal{B}(H)$ is also asymorphic to $\mathcal{B}(Z)$.\

\textit{Case 2:} $\gamma=\lambda^{+}$ for some cardinal $\lambda$.
We write $G$ as a union $G=\cup_{\alpha<\gamma}G_{\alpha}$ of an
increasing family of subgroups such that $|G_{\alpha}|=\lambda$,
$|G_{\alpha+1}:G_{\alpha}|=\lambda$ for every $\alpha<\gamma$, and
$G_{\beta}=\cup_{\alpha<\beta}G_{\alpha}$ for every limit ordinal
$\beta$. Put $\kappa_{\alpha}=\lambda$ for every $\alpha<\gamma$.
By Theorem 2, $\mathcal{B}(G)$ is asymorphic to $\mathcal{B}(Z)$,
where $Z$ is defined by the family of parameters
$\{\kappa_{\alpha}:\alpha<\gamma\}$. Since $H$ admits a filtration
with the same family of parameters, $\mathcal{B}(H)$ is also
asymorphic to $\mathcal{B}(Z)$. 

\
This completes the proof.

\

\
It should be mentioned that Theorem 3 does not hold for countable
groups. By [2, Theorem 10.6], there exists a family $\mathcal{F}$
of countable locally finite groups such that any two groups from
$\mathcal{F}$ are non-asymorphic and
$|\mathcal{F}|=2^{\aleph_{0}}$. \

We do not know if Corollary 2 and Theorem 3 are true for groups of
singular cardinalities.

\
 
\begin{description}
\item \LARGE \textbf{References} \normalsize \item{[1]}  A. N.
Dranishnikov, \textit{Asymptotic topology}, Russian Math Surveys
55 (2000), no. 6, 1085-1129. \item{[2]}  Protasov I., Banakh T.
\textit{Ball Structures and Colorings of Graphs and Groups},
Matematical Studies Monograph Series, 11. L'viv: VNTL Publishers, 2003. \item{[3]}  Ihor
Protasov and Michael Zarichnyi, \textit {General Asymptology}, Mathematical Studies
Monograph Series, 12. L'viv: VNTL Publishers, 2007. \item{[4]}  John Roe,
\textit{Lectures on Coarse Geometry}, University Lecture Series,
31. Providence, RI: AMS, 2003.
\end{description}

\

\

\begin{flushright}
I.V. Protasov ($protasov@unicyb.kiev.ua$) \qquad \qquad \qquad
\qquad \qquad \qquad \qquad Department of Cybernetics, Kyiv
National University, Volodimirska 64, Kiev 01033, UKRAINE
\\
\end{flushright}
\begin{flushright}
Anastasiia Tsvietkova ($tsvietkova@math.utk.edu$)\qquad \qquad
\qquad \qquad \qquad \qquad \qquad  Department of Mathematics,
University of Tennessee, 121 Ayres Hall, Knoxville, Tennessee
37996-1300, USA
\end{flushright}

\end{document}